\documentclass[11pt]{article}
\usepackage{amssymb}
\usepackage{amsthm}
\usepackage{amsmath}

\usepackage{epsfig}

\usepackage{fullpage}

\begin{document}

\title{\bf Limiting distribution of last passage percolation models
}

\author{{\bf Jinho Baik\footnote{Department of Mathematics,
University of Michigan, Ann Arbor, MI 48109 USA,  E-mail:
baik@umich.edu}}}



\maketitle

\begin{abstract}
We survey some results and applications of last percolation models
of which the limiting distribution can be evaluated.
\end{abstract}

\section{Introduction}

\begin{figure}[ht]
\centerline{\epsfxsize=5cm\epsfbox{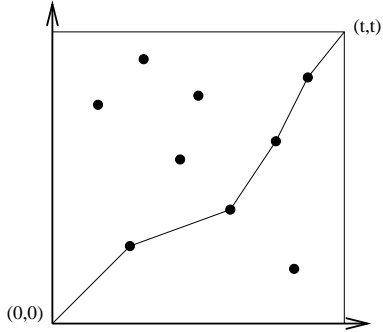}} \caption{Poisson
points and a longest up/right path from $(0,0)$ to $(t,t)$. In
this example, $L(t)=4$.}\label{fig:Poisson}
\end{figure}

Consider a Poisson process of rate $1$ in $\mathbb{R}_+^2$. If one
is only interested in the points in a square $(0,t)\times (0,t)$,
it is equivalent to think of picking $\mathcal{N}$ points at
random in the square where $\mathcal{N}$ itself is a random
variable such that
\begin{equation}
  \mathbb{P}\bigl(\mathcal{N}=N\bigr) = \frac{e^{-t^2}(t^2)^N}{N!}, \qquad
  N=0,1,2,\cdots.
\end{equation}
Given a realization of the process, an up/right path from $(0,0)$
to $(t,t)$ is a piecewise linear curve starting from $(0,0)$ to
$(t,t)$ joining Poisson points such that the slope, where defined,
is positive. The length of an up/right path is defined by the
number of points on it. Let $L(t)$ be the \textit{length of the
longest up/right path} from $(0,0)$ to $(t,t)$, making it a random
variable. See Fig.~\ref{fig:Poisson} for an example of a last
longest up/rigth path.
This is a type of last passage percolation model which is to find
a path that maximizes the total weight (passage time) in a random
environment.

Various statistics of $L(t)$ as $t\to\infty$ have been of
interest. The basic law is\cite{BDJ}
\begin{equation}\label{eq1}
  \lim_{t \to\infty} \mathbb{P}\biggl( \frac{L(t)-2t}{t^{1/3}} \le x\biggr)
  = F(x) := \exp \biggl( - \int_x^\infty (s-x) q^2(s) ds \biggr)
\end{equation}
where $q(x)$ is the (unique) solution\cite{HM} to the Painle\'e II
equation
\begin{equation}
  q''= 2q^3+xq
\end{equation}
satisfying the condition
\begin{equation}
  q(x) \sim -Ai(x) , \qquad x\to +\infty
\end{equation}
where $Ai$ denotes the Airy function. Convergence of all the
moments is also proved in the same paper. The limiting
distribution $F(x)$ is called the Tracy-Widom distribution, which
will be discussed further in section~\ref{sec:RM}.

This article intends to provide some motivations and applications
of the above Poisson last passage percolation model, and discuss
the result \eqref{eq1}.

\section{Motivations and Applications}


\subsection*{Random permutation}
Given a realization of the Poisson points in the square
$(0,t)\times (0,t)$, suppose one label them as $(x_i,
y_{\pi(i)})$, $i=1,\dots, \mathcal{N}$ such that
$x_1<x_2<\cdots<x_{\mathcal{N}}$. Note that with probability $1$
no two points have the same $x$- or $y$-coordinate. Then the
indices of the $y$-coordinates of the points generate a
permutation $\pi$. Moreover, an up/right path is mapped to an
increasing subsequence of the corresponding permutation. In the
example of Fig.~\ref{fig:Poisson}, the associated permutation is
$629473518$ and the increasing subsequence corresponding to the
indicated up/right path is $2358$. Therefore denoting the length
of the longest increasing subsequence of $S_N$ by $L_N$ and
recalling the property of the Poisson process discussed earlier,
we find that
\begin{equation}
  \mathbb{P}\bigl( L(t) = \ell) = \sum_{N=0}^\infty
  \frac{e^{-t^2}(t^2)^N}{N!} \mathbb{P}\bigl( L_N=\ell \bigr),
  \qquad \ell \in \mathbb{Z}.
\end{equation}
Using this formula and the so-called de-Poissonization
lemma\cite{Jo2} one can extract from the result \eqref{eq1}
that\cite{BDJ}
\begin{equation}\label{eq2}
  \lim_{N\to\infty} \mathbb{P} \biggl( \frac{L_N-
  2\sqrt{N}}{N^{1/6}} \le x \biggr) = F(x).
\end{equation}

The problem of finding various statistics of the longest
increasing subsequence of a random permutation in the large $N$
limit has been known as Ulam's problem since early 1960's. The
Poisson version of the model as in the last passage percolation
model above is sometimes called the Hammersley's process. See
e.g.~[1, 16] for more history of this combinatorial problem and
its applications.


\subsection*{Plancherel measure}
A partition of $N$ is a sequence of integers
$\lambda=\{\lambda_j\}_{j=1}^N$ such that $\lambda_1\ge \lambda_2
\ge\cdots\ge \lambda_N\ge 0$ and such that the sum of
$\lambda_j$'s is $N$. Given a partition $\lambda$ of $N$, let
$d_\lambda$ denote the number of standard Young tableaux of shape
$\lambda$ (see e.g.~[46] for definition). It is a basic result of
representation theory of the symmetric group that $d_\lambda$ is
the dimension of the irreducible representation of $S_N$
parameterized by $\lambda$. Hence the sum of $d_\lambda^2$ over
all partitions $\lambda$ of $N$ is equal to $N$. In view of this
identity, a natural probability on the set of partitions of $N$ is
the Plancherel measure defined by
\begin{equation}\label{eq:Plan}
  \mathbb{P} (\lambda) = \frac{d_\lambda^2}{N!}.
\end{equation}

Now the famous Robinson-Schensted\cite{Sc} algorithm states that
one can uniquely associate a pair of Young tableaux of a partition
$\lambda$ of $N$ to each permutation $\pi\in S_N$. Moreover
$L_N(\pi)$ is equal to the largest part $\lambda_1$ of the
corresponding partition. This implies that the distribution of the
largest part $\lambda_1$ of a random partition of $N$ taken
according to the Plancherel measure \eqref{eq:Plan} is precisely
equal to the distribution of $L_N$ of a random permutation. Hence
the result \eqref{eq2} yields the limiting law for $\lambda_1$.

The asymptotic statistics of other parts $\lambda_2,\lambda_3,
\dots$ have also been studied. See e.g.~[32, 5, 37, 14, 27, 6,
49].

\subsection*{Polynuclear growth (PNG) model}
\begin{figure}[ht]
\centerline{\epsfxsize=8.5cm\epsfbox{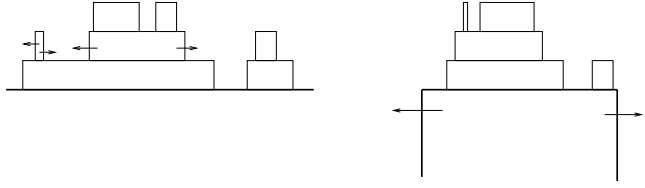}} \caption{A snapshot
of the height function of a PNG model. The picture on the right is
of droplet initial condition.}\label{fig:PNG}
\end{figure}

Consider a one-dimensional flat substrate. Suppose that there
occur random nucleation events, which is a Poisson process in the
space-time plane. If a nucleation occurs at $(x_0,t_0)$, an island
of height $1$ with zero width is created at position $x_0$ at time
$t_0$. As time increases, the island grows laterally in both
directions with speed $1$ while keeping its height. Often two
growing islands of same height collide. In that case they form one
island and the edges of the amalgamated island keep growing with
the same speed $1$. Note that nucleations can occur on top of an
existing island, and hence new islands may be created on an old
island. Let $h(x,t)$ be the height of PNG model at position $x$ at
time $t$. Thus $h(x,t)$ is a piecewise constant function in $x$
for fixed $t$. An example of the graph of $h(x,t)$ is in
Fig.~\ref{fig:PNG}.

One could impose various initial conditions. For now, suppose that
the nucleation events occur only for $|x|<t$. A different way of
describing this condition is that at time $t=0$, the substrate
satisfies $h(0,0)=0$ and $h(x,0)=-\infty$ for all $x\neq 0$, and
as time increases, the base substrate itself grows laterally of
speed $1$ and nucleations occur only on top of the base substrate
(and of course on top of islands on the base substrate). This
condition is called the \textit{droplet case}. See the picture on
the right in Fig.~\ref{fig:PNG} for an example.

An observation by Pr\"ahofer and Spohn\cite{SpohnP1} is that the
PNG model with the droplet initial condition can be mapped to the
above Poisson last passage percolation model. Imagine the
space-time plane in which the Poisson points corresponding to
nucleation events are marked by dots. Due to the droplet
condition, the dots are in the forward-light-cone $|x|<t$, and as
the growth speed of islands is $1$, the height at $(x_1, t_1)$
depends only on the nucleation events in the backward-light-cone
$|x-x_1|<t_1-t$. Therefore, $h(x_1, t_1)$ depends only on the
Poisson points in a rectangle of area $|t_1^2-x_1^2|/2$ . By
rotating the coordinates by $-\pi/4$, and by re-scaling, one ends
up with a Poisson process of rate $1$ in a square of area
$|t_1^2-x_1^2|/2$. Moreover, one can observe that the height
$h(x_1,t_1)$ is precisely equal to the length of the longest
up/righ path.

By re-interpreting the result \eqref{eq1} in terms of the PNG
model using this identification, Pr\"ahofer and
Spohn\cite{SpohnP1} found the following result: for fixed $|c|<1$,
\begin{equation}\label{eq3}
  \lim_{t\to\infty} \mathbb{P}
  \biggl( \frac{h(ct,t)-\sqrt{2(1-c^2)}t}{((1-c^2)/2)^{1/6}t^{1/3}} \le x
  \biggr)= F(x).
\end{equation}
Note that the super-diffusive scaling $t^{1/3}$ is due to the fact
that the height at a position strongly depends on the heights at
the neighboring positions. Actually it has been believed that the
exponent $1/3$ should be universal for a wide class of
one-dimensional random growth models as long as the spatial
correlation is not too weak (see e.g.~[13]). Such models are said
to be in the KPZ universality class. In the famous work\cite{KPZ},
Kardar, Parisi and Zhang introduced a nonlinear stochastic
differential equation for the surface height $h(x,t)$ as a
continuum model for this class of random growth models, and
renormalization group analysis has suggested that the scaling
exponent should be $1/3$ for one-dimension case. The PNG model now
plays the role of the unique growth model for which the $1/3$
exponent can be rigorously proved. Moreover, one can even
establish the limiting distribution as in \eqref{eq3}.

One might ask what happens if the initial condition is changed. It
turned out that while the scaling exponent is the same, the
limiting distribution may change. For instance, instead of the
droplet condition, consider the flat initial condition: the
initial substrate is $\mathbb{R}$ and nucleation events could
occur at any position on it. Then for any fixed position $x$, one
finds (see [10, 40])
\begin{equation}
  \lim_{t\to\infty} \mathbb{P}
  \biggl( \frac{h(x,t)-\sqrt{2}t}{2^{-1/6}t^{1/3}} \le x
  \biggr)= \exp \biggl( - \frac12 \int_x^\infty q(s) ds \biggr)\cdot
  F(x)^{1/2}.
\end{equation}
As a second example, consider a PNG model on a half line $\{ x:
x\ge 0\}$. In this case, one could imagine the situation such that
extra nucleation events occur \textit{at $x=0$}, which corresponds
to a 1-dimensional Poisson process in time at $x=0$. In other
words, there is excessive creation of islands at the origin. Then
depending on the rate $\alpha$ of the creation of islands at the
origin, the height function could have different properties.
Indeed if $\alpha$ is big, then the height at $x=0$ is dominated
by the creation of islands at the origin, while if $\alpha$ is
small, then it is likely that the islands created at the origin
have little effect for the height at $x=0$. Thus one expects a
transition of $h(0,t)$ in $\alpha$, which is actually proved in
[10, 42]. On the other hand, at $x=ct$ for a fixed $c>0$, 
the height is proved to have the fluctuation law given by $F(x)$ for 
the choices $\alpha=0, 1$ \cite{SasamotoI}. The authors of \cite{SasamotoI}
also computed height fluction for the transitional case for $\alpha=0,1$ 
when 
$x \sim t^{2/3}$ in terms of a Fredholm determinant. It is yet to be 
seen to obtain a Painlev\'e II type formula for this determinant. 
For further references in this direction, see e.g.~[3,
10, 11, 12, 41, 42, 43]. We note that all these different initial
conditions have combinatorial meanings on random permutations and
Plancherel measure.

\subsection*{Particle/anti-particle process}
A different description of the PNG model\cite{SpohnPscale} is to
regard the right edge of an island as a particle and the left edge
of an island as an anti-particle. Hence there are right-moving
particles and left-moving anti-particles on the real line.
Creation of island corresponds to creation of particle and
anti-particle pair, and the fact that two islands stick together
when they meet implies that upon colliding, particle and
anti-particle annihilate each other. In this dynamic picture, the
height function is now equal to the total number of the particles
and the anti-particles that have crossed the given position up to
the given time.

\subsection*{Random vicious walks}
The combinatorics of the longest increasing subsequence and the
Plancherel measure have an interpretation as non-intersecting
paths, which is sometimes called vicious walks\cite{Fisher}. See
e.g.~[21, 22, 23, 29]  for reference.


\section{Random matrix and universality}\label{sec:RM}

The Tracy-Widom distribution $F(x)$ in \eqref{eq1} also appears in
a totally different subject; random matrix theory. The main
interest in the random matrix theory is the limiting statistics of
the eigenvalues as the size of matrix tends to infinity. Random
matrix theory has very diverse applications in both mathematics
and physics from the spectrum of heavy nuclei to the zeros of
Riemann-zeta function (see e.g.~[35, 17, 31]).

Of special interest is the Gaussian unitary ensemble (GUE) which
is the set of $N\times N$ Hermitian matrices $H$ equipped with the
probability measure\cite{Mehta}
\begin{equation}\label{eq:GUE}
  \frac1{Z_N} e^{-N tr(H^2)} dH
\end{equation}
where $dH$ is the Lebesque measure and $Z_N$ is the normalization
constant. In 1994 Tracy and Widom\cite{TW1} considered the
limiting distribution of the largest eigenvalue $\xi_{\max}(N)$ of
$N\times N$ Hermitian matrix taken from GUE and found that
\begin{equation}\label{eq4}
  \lim_{N\to\infty} \mathbb{P} \bigl( (\xi_{\max}(N) -
  \sqrt{2}))\sqrt{2}N^{2/3} \le x \bigr)
  =F(x)
\end{equation}
where $F(x)$ is precisely the same function in \eqref{eq1}. In
other words, upon proper scaling, the largest eigenvalue of a
random Hermitian matrix taken from GUE and the length of the
longest up/right path in the Poisson last passage percolation
model have the same limiting law.

It should be mentioned that the analysis of obtaining \eqref{eq1}
and \eqref{eq4} are independent. Especially it is not found
whether there is a direct relation between $L(t)$ and
$\xi_{\max}(N)$ for finite $t$ and $N$. Only in the limit
$t\to\infty$ and $N\to\infty$, two seemingly different quantities
have the same limit after independent computations. A framework to
understand this might be central limit theorem. In the classical
central limit theorem, the sum of $n$ independent identically
distributed random variables converges, after proper scaling, to
the Gaussian distribution as $n\to\infty$, irrelevant to the
detail of the random variable. The results \eqref{eq1} and
\eqref{eq4} have the similar feature that the two different
`random variables' share the same limiting distribution. So one
may expect that the Poisson percolation model and GUE are two
instances of models to which a nonlinear version of central limit
theorem is applied.
It is however not known whether there is indeed such a natural
nonlinear version of the central limit theorem with proper general
setting that includes the Poisson percolation model and the GUE.
But there are some `patches' of universality results in random
matrix theory, and in the remaining of this section, we discuss
some of them.

\subsection*{Universality of random matrices}
A generalization of GUE is the set of Hermitian matrices with the
density
\begin{equation}\label{eq:GUEV}
  \frac{1}{Z_N} e^{-N tr V(H)} dH
\end{equation}
where $V$ is an analytic function with sufficient decay
properties. We remark that then the eigenvalues $\xi_1\ge \xi_2\ge
\cdots\ge \xi_N$ have the density
\begin{equation}\label{eq:density}
  \frac1{Z_n'} \prod_{1\le i<j\le N} |\xi_i-\xi_j|^2 \cdot
  \prod_{j=1}^N w(x_j)
\end{equation}
for some constant $Z_N'$ and $w(x)=e^{-NV(x)}$. This can be
thought of as a Coulomb gas of particles in $\mathbb{R}$ with
logarithmic repulsion and external potential $w$. For general $V$,
the local statistics of the eigenvalues in the middle of the
limiting density of states are found to be independent of $V$ (see
[39, 15, 18]). Moreover, by using the results of [18] and [34],
one can also deduce that the limiting distribution of the largest
eigenvalue of random Hermitian matrix taken according to the
probability \eqref{eq:GUEV} is generically given by the
Tracy-Widom distribution $F(x)$. The special case of $V(x) = x^4 +
t x^2$ was obtained by Bleher and Its\cite{BIts}.

The GUE has an alternative definition and a different
generalization. Namely, GUE is the set of Hermitian matrices with
independent, except for the Hermitian condition, complex Gaussin
entries. It is direct to check from this definition that the
density of the matrix is precisely \eqref{eq:GUE}. A natural
generalization is then the random Hermitian matrices of
independent entries which are not necessarily Gaussian. This is
called the Wigner matrix\cite{Mehta}. For Wigner matrix, the
largest eigenvalue is found to have the same limit $F(x)$
again\cite{Soshnikov} (see also [28] for a result regarding the
eigenvalues in the middle of the limiting density of states).
Hence the limiting law \eqref{eq4} still holds true for two
different generalizations of GUE, one of density function
\eqref{eq:GUEV} and the other of independent entries.

\subsection*{Growth models and discrete orthogonal polynomial ensembles}
In addition to the Poisson percolation model, there are a few more
isolated examples of percolation models for which the scaling
limit \eqref{eq1} can be obtained. Consider the lattice sites
$\mathbb{N}^2$. Suppose that we assign a random variable $X(i,j)$
at each site $(i,j)\in \mathbb{N}^2$. We further assume that
$X(i,j)$ are independent and identically distributed. Consider an
up/right path $\pi$ from the site $(1,1)$ to $(M,N)$, which is a
collection of neighboring sites $\{(i_k,j_k)\}$ such that
$(i_{k+1}, j_{k+1})-(i_k,j_k)$ is either $(1,0)$ or $(0,1)$. Let
$\Pi(M,N)$ denote the set of up/right paths from $(1,1)$ to
$(M,N)$, and define
\begin{equation}
  L(M,N):= \max_{\pi\in \Pi(M,N)} \biggl\{ \sum_{(i,j)\in \pi} X(i,j)
  \biggr\}.
\end{equation}
If $X(i,j)$ are positive random variables, an interpretation is
that $X(i,j)$ is the passage time at the site $(i,j)$, and
$L(M,N)$ is the last passage time to travel from $(1,1)$ to
$(M,N)$ along an admissible up/right path. The Poisson percolation
model is a continuum version of this more general directed last
passage percolation model.

For general random variables $X(i,j)$, the scaling limit law
\eqref{eq1} is an open problem, but when the random variable
$X(i,j)$ is either geometric or exponential, \eqref{eq1} is
proved\cite{kurtj:shape}. Also if the definition of the up/right
paths is modified, there are a few more isolated cases for which
\eqref{eq1} is obtained (see e.g.~[48, 27, 2]). However all the
cases such that \eqref{eq1} is proved share the common feature
that they all have an interpretation as a version of the longest
increasing subsequence of a random (generalized) permutation, and
all of them have the same algebraic structure (see e.g.~[38, 9]).
Hence it is an open question to prove the universality result
\eqref{eq1} for a general random variable $X(i,j)$ which does not
have such a structure.

The geometric percolation model above has an alternative
representation. Consider the density function on the set of
particles $\xi_1>\xi_2>\cdots>\xi_N\ge 0$, $\xi_j\in\mathbb{N}\cup
\{0\}$, given by \eqref{eq:density}. The only change is that the
`particles' $\xi_j$ lie on a discrete set instead of a continuous
set. Sometimes this is called the `discrete orthogonal polynomial
ensemble', while \eqref{eq:density} with continuous weight $w$ is
called the `continuous orthogonal polynomial ensemble'. A result
of Johansson\cite{kurtj:shape} is that for the geometric
percolation model, $L(M,N)$ has the same distribution (except for
a minor translation change) as the `largest particle' $\xi_1$ in
the discrete orthogonal polynomial ensemble with the special
choice $w(x)=\binom{x+M-N}{x} q^{x}$ (assuming that $M\ge N$).

Hence one may wonder whether the discrete orthogonal polynomial
ensemble with general weight $w$ have universal properties just
like its continuous weight counterpart \eqref{eq:density}. This is
indeed proved to be the case in [7, 8] for a wide class of
discrete weight $w$. The analysis extends the Deift-zhou
steepest-descent method\cite{DZ1} of Riemann-Hilbert analysis to
the discrete interpolation setting (see also [33, 36]). For
general discrete weights $w$, it is not clear if discrete
orthogonal polynomial ensembles have any percolation-type
interpretation, but the universality of \eqref{eq:density} for
both continuous and discrete weights provides a linkage of the
result \eqref{eq4} for GUE and the result \eqref{eq1} for the
geometric percolation model.

\section*{Acknowledgments}
The author would like to thank Percy Deift and Michio Jimbo for
kindly inviting him to the special session of integrable systems
in the International Congress on Mathematical Physics. The work is
supported in part by NSF Grant \# DMS-0208577.

%
%
%
%

\end{document}